\documentclass[10pt,english]{amsart}
\usepackage{leqno,latexsym,babel,amssymb}
\usepackage[latin1]{inputenc}
\usepackage{epsfig,graphicx}
\def\hf{\frac{1}{2}}
\def\th{\theta}
\def\tho{\theta_{\text{odd}}}

\def\cM{{\mathcal M}}
\newtheorem{thm}{Theorem}[]
\newtheorem{lem}{Lemma}[]
\newtheorem{prop}{Proposition}[]
\newtheorem{cor}{Corollary}[]
\def\ni{\noindent}
\def\byx{\hfill{$\square$}}

\def\vre{\varepsilon}

\begin{document}

\title{On the congruence $ax+by \equiv 1$ modulo $xy$}

\abstract We give bounds on the number of solutions to the
Diophantine equation $(X+1/x)(Y+1/y) = n$ as $n$ tends to
infinity. These bounds are related to the number of solutions to
congruences of the form $ax+by \equiv 1$ modulo $xy$.

\endabstract


\author{J. Brzezi\'nski}

\author{W. Holszty\'nski}
\address{400 E.Remington Dr\\
\#D233 Sunnyvale, CA 94087, USA } \email{wlod2@earthlink.net}

\author{P. Kurlberg}
\address{Department of Mathematics, Chalmers University of Technology
\hfill\break  and G\"oteborg University, S--41296 G\"oteborg,
Sweden} \email{jub@math.chalmers.se, kurlberg@math.chalmers.se}

\subjclass{Primary 11D45, Secondary 11A25, 11D72}

\keywords{Diophantine equation, linear congruence, divisor
function}

\maketitle


\section{INTRODUCTION}

Erik Ljungstrand has asked the first author about estimates of the
number of solutions to the equation
$$n=\left(X+\frac{1}{x}\right)\left(Y+\frac{1}{y}\right),\eqno(1)$$
where $n,X,x,Y,y$ are positive integers satisfying $n >1$, $x
> 1$ and $y
> 1$. His computations suggested that the number of such
solutions, when symmetric solutions obtained by transposing
$(X,x)$ and $(Y,y)$ are identified, is always less than $n$.

It is easy to see that $y$ divides $xX+1$ and $x$ divides $yY+1$.
Denoting the corresponding quotients by $b$ and $a$, we get the
following system:
$$ax = yY+1,$$
$$by =xX+1,$$

\ni
 where $ab = n$. Thus

$$ax \equiv 1 \pmod y \hskip2em \text{and} \hskip2em by \equiv 1 \pmod
x. \eqno(2)$$

It is clear that the integers $x,y$ satisfying these congruences
are relatively prime, and the system is equivalent to

$$ax + by \equiv 1 \pmod {xy}. \eqno(3)$$

It is also clear from the equations above that $x \neq y$, so when
counting the solutions, we may assume $x < y$. It is not difficult
to see that the problem of finding all solutions to equation (1)
with $1 < x < y$ is equivalent to the problem of finding all
solutions to the systems of linear congruences (2) for all $a,b$
such that $ab = n$ with $x,y$ satisfying the same conditions (see
Section 2).

One of the aims of the present paper is to prove E. Ljungstrand's
observation concerning the number $f(n)$ of solutions to equation
(1). The proof is a combination of an estimate of $f(n)$ (see
Theorem 3) proving the result for relatively big values of $n$ and
a portion of numerical computations, which together prove the
inequality $f(n) < n$ for all $n$. The systems of linear
congruences (2) or the congruence (3) (for fixed $a,b$) seem to be
interesting on their own rights. In the paper, we study the sets
of solutions to these congruences and give some estimates for
their size both from above and below. We give also a reasonably
effective algorithm for finding all solutions of (1) in positive
integers and attach some numerical results. In the last part of
the paper, we study the arithmetic mean of the function $f(n)$ and
give some lower and upper bounds for its size.


\section{CONGRUENCES}

Our objective is to estimate the number of solutions with $x,y
> 1$ to the congruence $ax + by \equiv 1 \pmod {xy}$ when $ab = n$ is fixed.

\begin{thm} Let $a,b$ be fixed positive integers and $ab = n > 1$. Let
$\rho(a,b)$ denote the number of pairs $(x,y)$ of integers $x,y$
such that $xy \mid ax + by -1$, $1 < x < y$. Then for every $n
\geq 1$ and for every real number $1 \leq \alpha \leq \sqrt{n}$,

$$\rho(a,b) < \frac{1}{\alpha}\sqrt{n}\log(n) +
2\left(1+\frac{0.6}{\alpha}\right)\sqrt{n} +
\frac{(2n-1)\alpha}{2\sqrt{n}-\alpha}.$$
\end{thm}

Before we prove the Theorem, we need two preparatory results. Let
$\theta(n)$ denote the number of divisors to $n$.

\begin{lem} Let $n \geq 22$ be a natural number and $1 \leq
\alpha \leq \sqrt{n}$ a real number. Then

$$\sum_{k=1}^{\frac{1}{\alpha}\root \of n} \theta(n-k) <  \frac{1}{\alpha}\sqrt{n}\log(n) +
2\left(1+\frac{0.6}{\alpha}\right)\sqrt{n}.$$
\end{lem}

\ni {\bf Proof.} We have (see e.g. [2], p. 347):

$$\sum_{k=1}^{\frac{1}{\alpha}\sqrt{n}}\theta(n-k) = \sum_{k=1}^{\frac{1}{\alpha}\sqrt{n}}
\sum_{d \mid n-k}1 \leq 2 \sum_{k=1}^{\frac{1}{\alpha}\sqrt{n}}
\sum_{\substack{d \mid n-k \\ 1 \leq d \leq \sqrt{n}}} 1 \leq 2
\sum_{d=1}^{\sqrt{n}} \left( \frac{\frac{1}{\alpha}\sqrt{n}}{d}+1
\right)$$

$$\leq \frac{2}{\alpha} \sqrt{n} (\log\sqrt{n} + 0.6) + 2\sqrt{n}
 = \frac{1}{\alpha}\sqrt{n} \log n + 2(1+\frac{0.6}{\alpha}) \sqrt{n},$$

\noindent where the last inequality follows noting that
$(\sum_{1}^{n}\frac{1}{k})-\log n$ is decreasing and less than 0.6
when $n \geq 22$.

\hfill \byx


\begin{lem} Let $a,b,x,y$ be positive integers such that $ab = n$, $ax
\equiv 1 \pmod y$, $by \equiv 1 \pmod x$ and $x,y > 1$. Let $ax -1
= yY$, $by-1 = xX$ and $ax+by-1 = kxy$. Then

\ni
 $($a$)$ $k = n - XY$,

\ni
 $($b$)$ $x = \frac{b+Y}{k}$ and $y = \frac{a+X}{k}$,

\ni $($c$)$ $\max (x,y) \leq \frac{2n-1}{2k-1}$,

\ni $($d$)$ $k \leq \frac{n+1}{3}$.

\end{lem}

\ni {\bf Proof.} We have

$$xyXY = (ax-1)(by-1) = abxy - ax - by +1 = abxy - kxy.$$

\ni Dividing by $xy$, we get (a). Now $ax-yY = by - xX$ gives
$x(a+X) = y(b+Y)$, so $\frac{a+X}{y} = \frac{b+Y}{x}$. But

$$kxy= ax+by - 1 = \left(\frac{ax-1}{y}+b\right)y = (Y+b)y$$

\ni shows that both fractions are equal to $k$, which proves (b).
We have

$$ky = a + X = a + \frac{ab-k}{Y} \leq ab + \frac{ab-k}{b+Y-1}
\leq ab + \frac{ab-k}{2k-1} = \frac{(2ab-1)k}{2k-1},$$

\ni where the last inequality follows from $b+Y = kx \geq 2k$, and
the first is equivalent to

$$ab -a = a(b-1) \geq \frac{ab-k}{Y} - \frac{ab-k}{b+Y-1} =
\frac{ab-k}{Y} \cdot \frac{b-1}{b+Y-1} = X \frac{b-1}{kx-1},$$

\ni that is, $a(kx-1) \geq X$, when  $b \neq 1$. This is
equivalent to $akx \geq a+X = ky$, which immediately follows from
$ax = yY+1
> y$. By symmetry, we get the corresponding
inequality with $y$ replaced by $x$, which proves (c).

Since $x,y \geq 2$ and, of course, $x \neq y$, we have $\max (x,y)
\geq 3$. Thus (c) implies (d).

\byx

\ni {\bf Proof of Theorem 1.}   Let $1 < x < y$ be integers such
that $xy \mid ax + by -1$. Notice that given $y$ there is only one
$x$ satisfying the necessary condition $ax \equiv 1 \pmod y$ and
therefore at most one pair $(x,y)$ such that $xy \mid ax + by -1$.

Using notations from Lemma 2, we have $XY = ab-k = n - k < n$.
Observe that $X$ and $Y$ are positive, since $x >1$ and $y > 1$.
We consider contributions to the numbers of solutions in two
cases.

 First of all, let $k \geq \frac{1}{\alpha}\sqrt{n}$,
 where $1 \leq \alpha \leq \sqrt{n}$.
  Then according to Lemma 2 (c), we get

$$y \leq \frac{2n-1}{2k-1} \leq
\frac{2n-1}{\frac{2}{\alpha}\sqrt{n}-1} =
\frac{(2n-1)\alpha}{2\sqrt{n}-\alpha}$$

\noindent Since every $y$ gives at most one $x$, we have less than
$\frac{(2n-1)\alpha}{2\sqrt{n}-\alpha}$ possibilities for $(x,y)$
in this case.

Assume now that $k < \frac{1}{\alpha}\sqrt{n}$ is fixed. Then,
since $X \mid n-k$, we get at most $\theta(n-k)$ possibilities for
its choice. But $k$ and $X$ uniquely define $y$, and consequently,
$x$. Therefore the number of possibilities for $(x,y)$ in this
case is at most $\sum_{k=1}^{\frac{1}{\alpha}\root \of n}
\theta(n-k)$, which according to Lemma 1 is less than:

$$\frac{1}{\alpha}\sqrt{n}\log(n) +
2\left(1+\frac{0.6}{\alpha}\right)\sqrt{n}.
$$

\noindent
 Thus the total number of possible $(x,y)$ is at most:

$$\frac{1}{\alpha}\sqrt{n}\log(n) +
2\left(1+\frac{0.6}{\alpha}\right)\sqrt{n} +
\frac{(2n-1)\alpha}{2\sqrt{n}-\alpha}.$$  \byx

Notice that if we fix $k < \frac{1}{\alpha}\sqrt{n}$ and choose
$X$ as a divisor to $n-k$, then $x$ and $y$ are uniquely
determined regardless of whether $x < y$ or $x > y$. In fact, $k$
and $X$ uniquely determine $y$, $Y$ (from $XY = n-k$) and,
consequently, $x$ from Lemma 2 (b). Thus if we are interested in
the total number of solutions to (3) without the assumption $x <
y$, then we have to count twice the number of solutions
corresponding to $k \geq \frac{1}{\alpha}\sqrt{n}$ (they may
correspond to $x < y$ or $x > y$) plus the number of solutions
corresponding to $k < \frac{1}{\alpha}\sqrt{n}$. Thus we have

\ni {\bf Theorem 1'.} {\it Let $a,b$ be fixed positive integers
and $ab = n$. Let $\rho'(a,b)$ denote the number of pairs $(x,y)$
of integers $x,y$ such that $xy \mid ax + by -1$, $x,y > 1$. Then
for every integer $n \geq 1$ and every real $1 \leq \alpha \leq
\sqrt{n}$,

$$\rho'(a,b) < \frac{1}{\alpha}\sqrt{n}\log(n) +
2\left(1+\frac{0.6}{\alpha}\right)\sqrt{n}+\frac{2(2n-1)\alpha}{2\sqrt{n}-\alpha}.$$
} \byx

For completeness of our discussion of the congruence (3), we note:

\begin{prop} The congruence $ax + by \equiv 1 \pmod {xy}$
has infinitely many solutions in positive integers $x,y$ if and
only if $a = 1$ or $b=1$.
\end{prop}

\ni  {\bf  Proof.} As we already know, there is only finitely many
solutions with $x,y > 1$. Therefore, if we have infinitely many
solutions, then in infinitely many of them $x = 1$ or $y = 1$. If
for example, $x=1$ then infinitely many $y$ divide $a-1$, so $a =
1$. The converse is trivial.\byx


\section{THE EQUATION}

In this section, we discuss the number of solutions to equation
(1), give an estimate of it and prove that for big values of $n$,
it is always less than $n$. First we note:

\begin{thm} $($a$)$ The solutions $(X,x,Y,y)$ to
the equation $(1)$ with $1 < x < y$
are in a one-to-one correspondence with the quadruples $(x,y,a,b)$
such that $ab = n$, $1 < x < y$ and $ax+by \equiv 1 \pmod {xy}$.

\ni $($b$)$ The solutions $(X,x,Y,y)$ to the equation $(1)$  with
fixed value $k = n - XY> 0$,  $x ,y
> 1$, $X \leq  Y$, and $x<y$ if $X = Y$, are in a one-to-one
correspondence with the set of the quadruples $(X,Y,a,b)$
satisfying

$$n = ab > n-k = XY, \hskip2em
k \mid \gcd(a+X,b+Y),\eqno(4)$$

\noindent where $a+X > k, b+Y > k$, $X \leq Y$ and $a < b$ if $X =
Y$. Moreover, for every solution $(X,x,Y,y)$ to the equation
$(1)$, $x = \frac{b+Y}{k}$ and $y = \frac{a+X}{k}$.

\end{thm}

\ni  {\bf  Proof.} (a) As noted in the introduction, a solution
$(X,x,Y,y)$ to equation (1) with $1 < x < y$ gives the congruence
$ax+by \equiv 1 \pmod {xy}$, where $a = \frac{yY+1}{x}$ and $b =
\frac{xX+1}{y}$, $ab = n$. Conversely, if $(x,y)$ is a solution to
$ax+by \equiv 1 \pmod {xy}$, where $ab = n$ and $1 < x < y$, then
we easily check that $(X,x,Y,y)$ with $X = \frac{by-1}{x}$ and $Y
= \frac{ax-1}{y}$ is a solution to equation (1).

(b) Let $(X,x,Y,y)$ be a solution to equation (1) with $k = n -
XY$, $x,y
> 1$, $X \leq Y$, and $x < y$ if $X = Y$.
Then with $a,b$ as above, we get a quadruple $(X,Y,a,b)$.
According to Lemma 2, $n = ab  > n-k = XY$, $k \mid
\gcd(a+X,b+Y)$, and $x = \frac{b+Y}{k}$, $y = \frac{a+X}{k}$.
Hence $x,y
> 1$ imply $a+X > k$ and $b+Y > k$. Moreover, if $X = Y$, then $x
< y$ gives $a < b$.

Conversely, if $(X,Y,a,b)$ is any quadruple satisfying the
conditions in (b), then we get $(X,x,Y,y)$, where $x =
\frac{b+Y}{k}$ and $y = \frac{a+X}{k}$, which is easily seen to be
a solution of the equation (1)  satisfying all the conditions in
(b). \byx

\medskip

\ni {\bf Remark 1.} Notice that the condition $k \mid
\gcd(a+X,b+Y)$ is equivalent to $\gcd(a+X,b+Y) = k$, since

$$k = ab - XY = (a+X)b - X(b+Y) $$

\noindent implies that $\gcd(a+X,b+Y) \mid k$. Moreover, if
$\gcd(n,k) = 1$, then the conditions $k \mid a+X$ and $k \mid b+Y$
are equivalent. In fact, $\gcd(n,k) = 1$ implies $\gcd(X,k) =
\gcd(b,k) = 1$, so the identity above implies the equivalence of
both conditions. Thus if $\gcd(n,k) = 1$, then in order to find a
solution to equation (1), it is sufficient to find factors $a$ of
$n$ and $X$ of $n-k$ such that $k \mid a+X$ with $a+X > k$ and
$\frac{n}{a}+\frac{n-k}{X} > k$. Then

$$\left(X,x =
\frac{\frac{n}{a}+\frac{n-k}{X}}{k},Y = \frac{n-k}{X},y =
\frac{a+X}{k} \right)$$

\ni is a solution. In particular, if $a = 1$, we obtain solutions
for every $k, X$ such that $\gcd(k,n) = 1$,

$$X \mid n-k,\,\,\,  k \mid X + 1 \,\,\,\,\text{and}\,\,\,\,
X+1 > k. \eqno(5)$$

\ni On the other hand, if $a = n$, we get solutions for $k,X$ such
that $\gcd(k,n) = 1$,

$$X \mid n-k,\,\,\, k \mid n+X\,\,\,\,\text{and}\,\,\,\,
1+\frac{n-k}{X} > 1. \eqno(6)$$

\ni  We shall use these observations frequently in Section 6.

Theorem 2 (a) implies that in order to estimate the number of
solutions to equation (1), we have to estimate the number $f(n) =
\sum_{ab =n}\rho(a,b)$ of solutions with $1 < x < y$ to all the
congruences $ax+by \equiv 1 \pmod {xy}$ when $ab =n$. It is well
known that for every $\vre
> 0$ there is a constant $C_{\vre}$ only depending on $\vre$ such
that $\theta(n) \leq C_{\vre} n^\vre$. Applying this fact and
Theorem 1, we get a bound on $f(n)$ depending on $n, \alpha$ and
$\vre$. However, we can get a somewhat sharper estimate noting
that we can only use one of the congruences $ax+by \equiv 1 \pmod
{xy}$ and $bx+ay \equiv 1 \pmod {xy}$, but instead, taking all
possible solutions with $x,y > 1$ (that is, removing the
assumption $x < y$). In fact, it is clear that $(x,y)$ solves the
first congruence if and only if $(y,x)$ solves the second one. In
such a way, we can use the estimate from Theorem 1', but only for
the pairs $a,b$ with $ab =n$ and $a \leq b$. The number of such
pairs is $\frac{1}{2}\theta(n)+\epsilon_n$, where $\epsilon_n = 0$
if $n$ is not a square and $\epsilon_n = \frac{1}{2}$,  when $n$
is a square. This gives the following result:

\begin{thm} Let $f(n)$ denote the number of solutions to the
equation $(1)$ and let

$$g(n,\alpha) = \frac{1}{\alpha}\sqrt{n}\log (n)
+2\left(1+\frac{0.6}{\alpha}\right)\sqrt{n} +
\frac{2(2n-1)\alpha}{2n -\alpha\sqrt{n}}.$$

\ni Then for every $\vre > 0$ and any real $1 \leq \alpha \leq
\sqrt{n}$ there is a constant $C_{\vre}$ such that

$$f(n) \leq \frac{1}{2}\theta(n)g(n,\alpha) \leq C_{\vre}
n^{\vre}\left(\frac{1}{2\alpha}\sqrt{n}\log(n) +
\left(1+\frac{0.6}{\alpha}\right)\sqrt{n}+
\frac{(2n-1)\alpha}{2\sqrt{n}-\alpha}\right),$$

\ni when $n$ is not a square, and

$$f(n) \leq \frac{1}{2}(\theta(n)+1)g(n,\alpha) \leq  (C_{\vre}
n^{\vre}+1)\left(\frac{1}{2\alpha}\sqrt{n}\log(n) +
\left(1+\frac{0.6}{\alpha}\right)\sqrt{n}+
\frac{(2n-1)\alpha}{2\sqrt{n}-\alpha}\right),$$

\ni when $n$ is a square. In particular, if $n$ is sufficiently
big then $f(n) < n$.
\end{thm}


\section{AN ALGORITHM}

We can now construct a reasonably efficient algorithm for
computing the number of solutions $(X,x,Y,y)$ to equation (1)
following their description in Theorem 2 (b).

First of all, write down the divisor list of $n$. For each divisor
$a$ of $n$ and for all integers $X$ such that $1 \leq X <
\sqrt{n}$, repeat the following: Compute all the divisors $k$ of
$a+X$, for each $k$, check whether $Y = \frac{n-k}{X}$ and $x =
\frac{b+Y}{k}$, where $b = \frac{n}{a}$, are integers or not, put
$y = \frac{a+X}{k}$, $x = \frac{b+Y}{k}$ in the former case. If $X
= Y$ and $x > y$ replace $(x,y)$ by $(y,x)$.  Check whether $x
> 1$, $y
> 1$ and accept the quadruple $(X,x,Y,y)$ as a solution if all these
conditions are satisfied.

Theorem 2 (b) easily implies that this algoritm gives all the
solutions to equation (1) and every solution exactly once.

We are now ready for the numerical computations proving that the
number $f(n)$ of solutions to equation (1) is always less than
$n$.

As we noted before, for each $\vre > 0$ there is a constant
$C_\vre$ only depending on $\vre$ such that $\theta(n) \leq C_\vre
n^\vre$ for all $n \geq 1$. For simplicity, let $\vre =
\frac{1}{4}$ and denote by $C^\ast$ the least constant
corresponding to this value of $\vre$. It is easy to show that on
the positive integers the quotient

$$C(n) = \frac{\theta(n)}{n^{\frac{1}{4}}}$$

\ni attains its maximum value for $n = 21621600$, which gives
$C^\ast < C_0 =8.44697$.

According to Theorem $3$, if $n$ is not a square, we want to
decide when

$$f(n) \leq \frac{1}{2}\theta(n)g(n,\alpha) \leq \frac{1}{2\alpha}C^\ast n^{\frac{3}{4}}\log(n) +
\left(1+\frac{0.6}{\alpha}\right)C^\ast n^{\frac{3}{4}}
+\frac{(2n-1)\alpha}{2\sqrt{n}-\alpha}C^\ast n^{\frac{1}{4}} <
n.$$

Let

$$h(n,\alpha,C) = n^{\frac{1}{4}} - \frac{1}{2\alpha}C\log (n)
-\left(1+\frac{0.6}{\alpha}\right)C - \frac{(2n-1)\alpha}{2n
-\alpha\sqrt{n}}C.$$


Choose $\alpha = 2.95$. Then it is easy to check that
$h(n,\alpha,C^\ast) > h(n,\alpha,C_0) > 0$ when $n \geq 11 621
000$. By the definition of $C^\ast$, this shows that $f(n) < n$
for all $n \geq 11 621 000$ and it remains to check this
inequality for all $n < 11 621 000$.  In order to carry out the
numerical computation, we find all the numbers $n$ for which
$\frac{1}{2}\theta(n)g(n,\alpha) \geq n$. This happens when
$\theta(n)$ is ''big'', which occurs for $n$ having many small
prime factors. The computations give 6523 numbers in the interval
$[2\cdot 10^4,11621000]$: 3030 in $[2\cdot10^4,10^5]$, 3482 in
$[10^5,5\cdot10^6]$ and 11 in $[5\cdot 10^6, 11621000]$. The
numbers in the last interval are 5045040 (4559), 5266800 (4051),
5405400 (5069), 5569200 (4494), 5654880 (4534), 5765760 (5286),
6126120 (5211), 6320160 (5407), 6486480 (4333), 7207200 (6309),
8648640 (5330), where the number in the parenthesis is the
corresponding value of $f(n)$.

If $n$ is a square, then we repeat the same procedure as above
taking into account the extra term on the right hand side in the
second inequality in Theorem $3$. The bound 11621000 works in this
case as well, so we have to consider all squares less than this
bound (3408 numbers). Short computations show that there are 118
such squares for which the expression in the second inequality in
Theorem $3$ is not less than $n$ (the biggest one 1587600). For
these 118 numbers, we check by computer calculations that $f(n) <
n$.


\section{REDUCED SOLUTIONS}

The main aim of this section is a non-computational proof of the
inequality  $f(n) < n$ for the case when $n = p$ is a prime
number. We also give some estimates of the number $k = n - XY$ for
the solutions $X,x,Y,y$ to equation (1).

Let $X,x,Y,y$ be a solution to equation (1), which in this section
will be denoted by $n = [X,x,Y,y]$. Recall that $a,b$ denote
integers such that $ax = yY+1$ and $by = xX+1$. We say that a
solution $X,x,Y,y$ is {\it reduced} if $X < y$ and $Y < x$. The
reduced solutions are characterized in the following way:

\begin{prop}  Let $n = [X,x,Y,y]$. Then $X,x,Y,y$ is
reduced if and only if $XY = n-1$.

\end{prop}

\ni {\bf Proof.} If $X < y$ and $Y < x$, then Lemma 2 gives $kxy =
ax + by -1 = xX + Yy+1 < x(y-1)+y(x-1)+1 = 2xy+ 1-x -y < 2xy$.
Thus $k = n - XY = 1$. Conversely, if $XY = n-1$, then by Lemma 2,
$k = n - XY = 1$. This implies $X < y$ and $Y < x$, since
otherwise, $kxy = xX+yY+1 > xy$, that is, $k
> 1$.\byx

\begin{cor} The number of reduced solutions to the equation
$(1)$ is  $\frac{1}{2}\theta(n)\theta(n-1)$.

\end{cor}

\ni {\bf Proof.} If $X,y,Y,y$ is a reduced solution, then $ab =
n$, $XY = n-1$ and $k = 1$ according to Lemma 2. Thus each pair of
divisors to $n$ and $n-1$ defines a solution and every solution
gives such a pair of divisors. Of course, we have to divide by 2
the total number of such pairs in order to obtain each
desymmetrized solution exactly once. \byx

\begin{prop} If $p$ is a prime, then $f(p) < p$.
\end{prop}

\ni {\bf Proof.} According to Corollary 1, the number of reduced
solutions to $p = [X,x,Y,y]$ equals $\theta(p-1)$. Assume that the
solution $X,x,Y,y$ is not reduced. Without loss of generality,  we
may assume that

$$xX+1 = py \hskip2em \text{and} \hskip2em yY+1 = x.$$

\ni The second equation gives $Y < yY+1 = x$. Since the solution
is not reduced, we have $X > y$ (the equality is of course
impossible by the first equation). The second equation gives
$y|x-1$, so $y<x$. We also have $x < p$, since otherwise $py =
xX+1 > pX$ gives a contradiction. Thus $x$ belongs to the set
$\{3,\ldots,p-1\}$ with $p-3$ elements. Moreover, $py \equiv 1
\pmod x$ and $y < x$, so the congruence allows at most one $y$
giving a solution to the equation. If now $p \equiv 1 \pmod x$,
then $y$ must be equal to 1, which is impossible. Thus $x > 2$ can
not assume values dividing $p-1$. The number of such $x$ is
$\theta(p-1)-2$. Thus $x$ assumes at most

$$(p-3) - (\theta(p-1)-2) = p - \theta(p-1) - 1$$

\ni different values which give non-reduced solutions. According
to Corollary 2,  the number of reduced solutions is $\theta(p-1)$
so the total number of solutions is at most $p-1$. \byx


Every solution $X,x,Y,y$ to equation (1) has the corresponding
value of $k = n - XY$. By Proposition 3, $k=1$ corresponds to the
reduced solutions. For these solutions, $X$ and $Y$ must be the
least positive solutions to the congruences $xX \equiv -1 \pmod y$
and $yY \equiv -1 \pmod x$ when $x,y$ are fixed. All other
positive solutions to these congruences, with $x,y$ fixed, are
given by $X+ry$, $Y+sx$ where $r,s \geq 0$. Thus starting from $n
= [X,x,Y,y]$ with a fixed pair $x,y$, we get

$$N = [X+ry,x,Y+sx,y],$$

\ni where $N = (rx+b)(sy+a)$. The number $n = ab$ is the least
number for which such a (reduced) solution with fixed $x,y$
exists. We have $N-(X+ry)(Y+sx) = k + r +s$. In particular, if
$r=1, s=0$ or $r=0, s=1$, we get quadruples for which the
corresponding parameter $k$ decreases by 1:

$$[X,x,Y,y] \mapsto [X+y,x,Y,y], \qquad[X,x,Y,y] \mapsto [X,x,Y+x,y]. \eqno(7)$$

\ni We shall say that these two transformations  are elementary.
Thus we can describe the solutions for a given $n$ in the
following way:

\begin{prop} Every solution to $n = [X,x,Y,y]$ with $k =
n-XY
> 1$ can be obtained from a reduced solution to $m =
[X_0,x,Y_0,y]$ for some $m < n$, by successive use of $k-1$
elementary transformations $(7)$.

\end{prop}

\ni {\bf Proof.} If we have a solution $n = [X,x,Y,y]$ with $k =
n- XY$ and $k
> 1$, then the solution is not reduced, which means that $X > y$
or $Y > x$, since Lemma 2 implies immediately that the equalities
are impossible. If $X > y$, then we get $n-(yY+1) = [X-y,x,Y,y]$,
while $Y
> x$ gives $n - (xX+1) = [X,x,Y-x,y]$ both with the corresponding value of
$k' = [n-(xX+1)]-X(Y-x)= k -1$. This ``reduction process"{}
eventually leads to a reduced solution for a natural $m < n$ and
the same $x,y$. Starting from such a reduced solution and
reversing the process, we get the given solution $n = [X,x,Y,y]$
after $k-1$ steps. \byx

By Lemma 2 (d), $k \leq \frac{n+1}{3}$. Observe, that for $t\ge 1$
and $n=3t-1$, we have $n = [1,2,2t-1,3]$ and in this case, $k=n-XY
= \frac{n+1}{3}$.















\section{SOME ESTIMATES}

We wish to give upper and lower bounds on the number of solutions $f(n)$ to
equation (1) when $n$ is averaged over some interval.
 For simplicity, if $g,h$ are positive functions, we
write $g(n) \ll h(n)$ if there is a positive constant $C$
 such that $g(n) \leq Ch(n)$ for all sufficiently big natural $n$.

\begin{thm} There exist positive constants $C_1, C_2$ such that
for $T \geq 2$,

$$C_1 < \frac{\sum_{1}^{T}f(n)}{T\log^3 T} < C_2.$$

\end{thm}

In the proof we need the following result:

\begin{lem}
$ \sum_{n \leq T} \theta(n)\theta(n-1) = O(T \log^2 T).$
\end{lem}
\begin{proof}
If $m \leq T$, we have
$$
\theta(m) \leq 2 \sum_{\substack{l|m \\ l \leq  \sqrt{T}}} 1
$$
and thus
$$
\sum_{n \leq T} \theta(n)\theta(n-1) \leq 4 \sum_{n \leq T}
\sum_{\substack{ l|n \\ l \leq \sqrt{T}}} \sum_{\substack{ k|(n-1)
\\ k \leq \sqrt{T}}} 1
$$
$$
= 4 \sum_{ k,l \leq \sqrt{T}} |\{n \leq T : l|n, \,  k|(n-1)  \}|
= 4 \sum_{ k,l \leq \sqrt{T}} |\{d \leq T/l : k|(dl-1)  \}|
$$
$$
= 4 \sum_{ k,l \leq \sqrt{T}} |\{d \leq T/l : d \equiv l^{-1} \mod
k  \}| \leq 4 \sum_{ k,l \leq \sqrt{T}} \left( \frac{T}{kl} + 1
\right)
$$
$$
= O(T \log^2 T) + O(T) = O(T \log^2 T).
$$
\end{proof}

\ni {\bf Proof of Theorem 4.} With  the notations from the
introduction, given $x,y$,  let us choose $X_0$ and $Y_0$ such
that $xX_0 \equiv -1 \pmod y$, $0<X_0<y$, $yY_0 \equiv -1 \pmod x$
and $0<Y_0<x$. We want to count the number of integers $X,Y \geq
1$ such that $X \equiv X_0 \pmod y$, $Y \equiv Y_0 \pmod x$ and
$$\left(X+\frac{1}{x}\right)\left(Y+\frac{1}{y}\right)
= n \leq T,
$$
\ni when $x,y > 1$ are fixed. Noting that

 $$
XY < \frac{(Xx+1)(Yy+1)}{xy} < 4XY
$$
we will obtain lower bounds by estimating from below the number of
$X,Y$ such that $4XY \leq T$.

\ni The congruences $X \equiv X_0 \pmod y$, $Y \equiv Y_0 \pmod x$
are equivalent to $X,Y$ being of the form
$$
X = X_0 + ry, \quad Y = Y_0 + sx
$$
for $r,s$ non-negative integers. Thus it is enough to estimate
$$
\#\{ r,s\geq 0 : (X_0 + ry)(Y_0 + sx) \leq T/4 \},
$$
which, since $X_0<y$ and $Y_0 <x$, we may bound from below by
$$
\#\{r,s \geq 0 : (r+1)(s+1)xy \leq T/4 \}.
$$
This, in turn, is greater than
$$
\#\{r,s \geq 0: rs \leq \frac{T}{16xy} \} \sim \frac{T}{16xy}\log
\frac{T}{16xy}.
$$
Summing over $x,y \leq T^{1/3}$, we then find that there are
$$
\gg \sum_{x,y\leq T^{1/3}}  \frac{T}{16xy}\log \frac{T}{16xy} \gg
\log T \sum_{x,y\leq T^{1/3}}  \frac{T }{16xy}
$$
$$
\gg T \log T \sum_{x,y\leq T^{1/3}} \frac{1}{xy} \gg T \log^3 T
$$
ways of finding $x,y,X,Y$ such that
$$
n = \frac{(Xx+1)(Yy+1)}{xy} \leq T.
$$
In other words, on average, there are at least $C_1\log^3 T$
solutions for some $C_1 > 0$.

In order to prove the existence of an upper bound, we note first
that if $r = s = 0$, then the solution $(X_0,Y_0,x,y)$ is reduced.
For $n \leq T$ the number of reduced solutions is according to
Corollary 1 and Lemma 3,

$$\sum_{n \leq T}\frac{1}{2}\theta(n)\theta(n-1) = O(T\log^2T).$$

Assume now that $r \geq 1$ and $s = 0$. Then the number of
solutions $(X,Y,x,y)$ to equation (1) such that $n \leq T$ is less
than

$$
\#\{r,x,y >0 : XY = (X_0 + r y)Y_0 \leq T \} \leq \#\{r,x,y >0 : r
yY_0 \leq T \}.$$

\ni Since $r \geq 1$, we have $M:=  yY_0 \leq T$ and since $yY_0
\equiv -1 \pmod x$, we have

$$\#\{r,x,y >0 : r yY_0 \leq T \} \leq \sum_{M \leq T}\#\{r,x,y >0 :
x \mid M+1,\, y \mid M,\,  rM \leq T\}$$

$$\leq \sum_{M \leq T}\theta(M)\theta(M+1)T/M,$$

\ni which, by partial summation and Lemma 3, is $O(T \log^3T)$.

The case $r = 0, s  > 0$  follows in a similar way to the previous
one.

Finally, if $r,s > 0$ and $(X,Y,x,y)$ is a solution to (1) such
that $n \leq T$, then since $XY < T$, we get

$$
\#\{r,s,x,y >0 : XY = (X_0 + r y)( Y_0 + sx) \leq T \}
$$
$$
\leq \#\{r,s,x,y >0 : r y s x   \leq T \} = O(T \log^3 T).$$

In other words, on average, there are at most $C_2\log^3 T$
solutions for some $C_2 > 0$. \byx

What else can be said about the size of $f(n)$? For instance, how
close is $f(n)$ to its average? As the following figure shows,
$f(n)$ oscillates rather widely.

\begin{figure}[h]
\centering

\hspace{-3.8cm}\includegraphics[totalheight=7cm,
angle=-90]{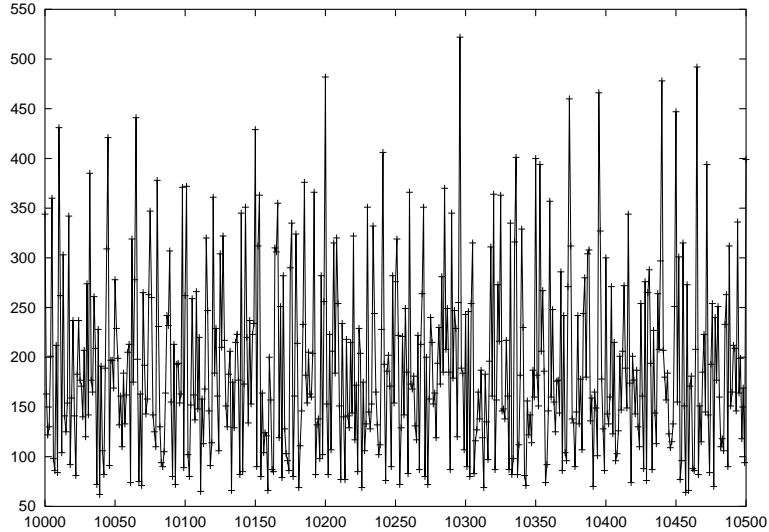}

\vspace{-2.7cm}

\hspace{3.5cm}
 \caption{The numbers of solutions to (1) for $10000
\leq n \leq 10500$}

\end{figure}





\ni Since there is $\frac{1}{2}\theta(n)\theta(n-1)$ reduced
solutions (see Corollary 1), it is clear that the order of
magnitude of $f(n)$ sometimes is larger than any power of $\log
n$. Moreover, there are other sources of large oscillations.

Let ${\mathcal M}(n,k)$ denote the number of solutions $X,x,Y,y$
to equation (1) such that $n - XY = k$. Of course,

$$f(n) = \sum_{k} \cM(n,k). \eqno(8)$$

\ni Taking into account the contribution to $f(n)$ from the number
of solutions with $k = 1$ and a similar contribution for $k = 2$
(see below Lemma 4), one might expect that the most significant
fluctuations of $f(n)$ depend on $\cM(n,k)$ for small values of
$k$. However, this is not the case as shown by the following
construction (we thank Andrew Granville for pointing this out to
us):  Fix an arbitrary $k$ and let $M > k$ be a large integer.
Choose $n = k + \prod p_i$, where $p_i$ are all primes such that
$p_i \equiv -1 \pmod k$ and $p_i \leq M$. Denote the number of
such primes $p_i$ by $\pi(M,k,-1)$. By the prime number theorem
for arithmetic progressions (see [1], Chap. 20 and 22):

$$\frac{c_kM}{\phi(k)\log M} \leq \pi(M,k,-1) \leq \frac{C_kM}{\phi(k)\log
M}$$

\ni for suitable positive constants $c_k,C_k$ only depending on
$k$. Now, half of the divisors to $n-k$ are congruent to $-1$
modulo $k$, so taking into account $(5)$, we get

$$f(n) \geq 2^{\pi(M,k,-1)-1}.$$

\ni Hence $\log f(n) \gg \pi(M,k,-1) \gg \frac{M}{\phi(k)\log M}$.
On the other hand,  since $\prod p_i \leq M^{\pi(M,k,-1)}$, we get

$$\log n \ll \sum \log p_i \ll
\frac{M}{\phi (k)},$$

\ni and similarly, $\log n \gg \frac{M}{\phi(k)}$, which implies
$\log M \ll \log\log n$. Thus

$$\log f(n) \gg \frac{M}{\phi(k)\log M} \gg \frac{\log n}{\log M}
\gg \frac{\log n}{\log\log n}.$$

\ni Hence

$$f(n) \gg \exp\left(\frac{c\log n}{\log\log n}\right)$$

\ni for some constant $c > 0$ only depending on $k$.

This shows that arbitrary $k$ may give ''big'' contribution to
$f(n)$ for a suitable $n$. It is also possible to show that the
contribution to $f(n)$ may come from many different values of $k$.
If $n+1$ has many different divisors, then according to $(6)$,
where we choose $X = 1$, each such divisor $k$ gives a solution to
the equation $(1)$.  Unfortunately, we are unable to prove that
$f(n) \rightarrow \infty$ when $n \rightarrow \infty$. What we
prove with ''some effort'' is

\begin{prop} If $n \geq 9$, then $f(n) \geq 8$.
\end{prop}

Let $\theta_{\text{odd}}(n)$ denote the number of odd divisors of
$n$. Then for $k = 2$, we have the following result:

\begin{lem} For $n\ge3,$ we have
$$\cM(n,2)= \begin{cases} \hf\th(n)\,\th(n-2)-1,& \text{if}\ n\ \text{is\
odd,}\\
\tho(n)\,\tho(n-2)-1,&\text{if}\ n\ \text{is\
even.}\\
\end{cases}$$
\end{lem}

\ni {\bf Proof.} In fact, if $n$ is odd, then according to Remark
1, we get all solutions to (1) taking any divisor $a$ to $n$ ($b =
\frac{n}{a}$) and any divisor $X$ to $n-2$ ($Y = \frac{n-2}{X}$)
such that $a+X
> 2$ and $b + Y> 2$. The number of pairs of
such divisors giving different quadruples $(X,x,Y,y)$ with $x < y
$ is $\hf\th(n)\,\th(n-2)$ and the only case when $a+X = 2$ or $b
+ Y = 2$ corresponds to the choice of $a = X = 1$ or $a = n$, $X =
n-2$, which gives only one quadruple with $x < y$. This proves the
first case.

If $n$ is even, let $n = 2^rm$, where $m$ is odd. One of the
numbers $n, n-2$ must be divisible by 4, so let us assume that $r
\geq 2$ (the case with $n-2$ divisible by 4 is considered in
similar way with the roles of $n, n-2$ interchanged). Thus $n-2 =
2(2^{r-1}m-1)$, and $\frac{n-2}{2}$ is odd. If $n-2 = XY$, then
exactly one of the factors $X,Y$ is even and the other one is odd.
Since $a+X$ and $b+Y$ are even, exactly one of the factors $a,b$
of $n = ab$ must be odd. Thus all the possibilities for the sums
$a+X$ and $b+Y$ are given by all the choices of the odd factors of
$n$ and $n-2$. Only one such choice gives $a+X=2$ or $b+Y =2$.
This proves the second case. \byx

Now we prove that

$$\text{if} \hskip1em
n > 11, \hskip1em {\text {then}} \hskip1em \cM(n,1) + \cM(n,2)
\geq 7. \eqno(9)$$

\ni First let $n$ be odd. Then

$$\cM(n,1) + \cM(n,2) =
\frac{1}{2}\theta(n)(\theta(n-1)+\theta(n-2))-1$$

\ni Since $n-1
> 4$ is even, $\theta(n-1) \geq 4$. Assume that $\theta(n) =
2$. Then $n$ is a prime. If also $\theta(n-2) = 2$, then $6 \mid
n-1$. Since  $n-1
> 6$, we have $\theta(n-1) \geq 6$, so $\cM(n,1) + \cM(n,2) \geq
7$. Assume now that $\theta(n-2) = 3$, that is, $n-2 = p^2$, where
$p > 3$ is a prime. Then  $3 \mid p^2 + 2 = n$, which is
impossible. Thus $\theta(n-2) \geq 4$, which gives $\cM(n,1) +
\cM(n,2) \geq 7$. Notice that if $n$ is a prime, $n-1$ twice a
prime and $n-2$ is a product of two different primes, then
$\cM(n,1) + \cM(n,2) = 7$. By Schinzel's conjecture (see [3]),
this situation happens for infinitely many $n$. If $\theta(n) >
2$, then it is easy to check that $\cM(n,1) + \cM(n,2) \geq 8$.

Assume now that $n$ is even, so

$$\cM(n,1) + \cM(n,2) = \frac{1}{2}\theta(n)\theta(n-1) +
\theta_{\text{odd}}(n)\theta_{\text{odd}}(n-2) -1.$$

We have $\theta(n) > 3$, since $n > 4$. Assume $\theta(n) = 4$.
Since $n
> 8$, we have $n = 2p$, where $p$ is an odd prime. If $n-1$ is a
prime, then $3 \mid n-2 = 2(p-1)$, so
$\theta_{\text{odd}}(n)\theta_{\text{odd}}(n-2) \geq 4$ and
$\cM(n,1) + \cM(n,2) \geq 7$. If $\theta(n) = 5$, then $n = 16$
and the claim follows by a direct computation. If $\theta(n) = 6$,
then $n = 32$ or $n = p^2q$ for two different primes $p,q$. If $p
=2$, then $n-2 = 2(2q-1)$ has at least two odd factors, so
$\cM(n,1) + \cM(n,2) \geq 7$. If $q = 2$ and $p = 3$, we check the
claim directly, and when $p > 3$, then $n-2 = 2(p^2-1)$ is
divisible by 3, so $\theta_{\text{odd}}(n)\theta_{\text{odd}}(n-2)
\geq 6$. If finally, $\theta(n) \geq 7$, then of course, the
inequality holds.

Now we prove that

$$\text{if} \hskip1em n > 12, \hskip1em \text{then} \hskip1em \cM(n,3)
\geq 1. \eqno(10)$$

Assume first that $3 \nmid n$ (so $3 \nmid n-3$) and let $n$ be
even. Then $n = 2m$ and $n-3 = 2m-3$. If for a prime $p \equiv 1$
(mod 3), $p \mid n-3$, then $p \geq 7$, so $X = p$, $Y =
\frac{n-3}{p}$, $x = \frac{2+Y}{3}
> 1$ and $y = \frac{m + X}{3} > 1$ (see Remark 1) give
a solution to equation (1). If for a prime $p \equiv 2$ (mod 3),
$p \mid n-3$, then $p \geq 5$, so $X = 1$, $Y = \frac{n-3}{p}$, $x
= \frac{p+Y}{3} > 1$ and $y = \frac{n+X}{3} > 1$ give such a
solution. If $n$ is odd, then $n-3$ is even and we repeat the same
arguments looking instead at the prime factors $p$ of $n$.

Let now $3 \mid n$. Let $n$ be even. Then $n = 3^s 2m$, where $3
\nmid m$, and $n-3 = 3(3^{s-1}2m-1) = 3r$. If $r$ has a prime
divisor $p \equiv 1$ or 2 (mod 3), we proceed exactly as in the
previous case above when $3 \nmid n$. Otherwise,  $2m-1$ is a
power of 3, so $n-3 = 3^{s+1}$ and $n = 3(3^s + 1)$. In this case,
$n$ must have a prime factor $p > 2$ congruent to 2 modulo 3 and
we get a solution to equation (1) as before.

If $n$ is odd, then $n-3$ is even and divisible by 3, so the
considerations are similar with the role of $n$ and $n-3$
interchanged.

Now the proof of Proposition 7 follows immediately from (8), (9),
(10) and by direct inspection of the cases $n = 9,10,11,12$. Still
more elaborate arguments show that $f(n) \geq 12$ if $n \geq 20$
(we thank Jerzy Browkin for sending us his proof of this result
and, in particular, for the proof of Lemma 4). \byx

\ni {\bf Remark 2.} It is no longer true that  $\cM(n,4) \geq 1$
for all sufficiently large $n$. If all primes dividing both $n$
and $n-4$ are congruent to 1 modulo 4, then by Remark 1, there are
no solutions to equation (1) with $k = 4$. In fact, this happens
for infinitely many $n$ by the following argument, for which we
thank Mariusz Ska{\l}ba. Let $m$ be a natural number such that
$m\not\equiv1\pmod{3}$ and put $n=2m^2+2m+5.$ Then

  $$n=(m-1)^2+(m+2)^2 \,\,\,\,\mbox{and}\,\,\,\, n-4=m^2+(m+1)^2$$

\ni are only divisible by primes congruent to 1 modulo 4.
\medskip

\ni {\bf Acknowledgments.} We thank Jerzy Browkin for many
valuable suggestions and helpful comments on an earlier version of
this paper. The third Author is partially supported by the Royal
Swedish Academy of Sciences and the Swedish Research Council.

\bigskip

\newpage

{\bf REFERENCES}

\bigskip

\begin{itemize}

\item[{[1]}] H. Davenport, {\it Multiplicative Number Theory}, Third edition,
Springer Verlag, 2000.

\item[{[2]}] G.H. Hardy, E.M. Wright, {\it An Introduction to the Theory of
Numbers}, Fifth edition, Oxford Science Pubications,  1979.

\item[{[3]}] A. Schinzel, W. Sierpi\'nski, {\it Sur certaines
hypoth\`eses concernant les nombres premiers}, Acta Arithmetica
4(1958), 185 -- 208.

\end{itemize}

\end{document}